\documentclass[11pt,a4paper]{amsart}
\usepackage{amssymb,amsmath,epsfig,graphics,mathrsfs}
\usepackage{amssymb,amsmath,epsfig,mathrsfs, enumerate, xparse, mathtools}
\usepackage[pagewise]{lineno}
\usepackage[pagebackref,colorlinks=true,linkcolor=blue,citecolor=blue]{hyperref}
\usepackage{graphicx}

\usepackage{fancyhdr}
\hypersetup{
 colorlinks   = true,
 urlcolor     = blue,
 linkcolor    = blue,
 citecolor   = red ,
 bookmarksopen=true
}

\usepackage{amsmath}
\usepackage{amsfonts}
\usepackage{amssymb}
\usepackage{amsthm}
\usepackage{epsfig,graphics,mathrsfs}
\usepackage{graphicx}

\usepackage{latexsym} 
\usepackage{amssymb,amsmath}
\usepackage{amsthm}
\usepackage{longtable,booktabs,setspace} 
\usepackage{url}

\usepackage[usenames, dvipsnames]{color}

\newcommand*\MSC[1][1991]{\par\leavevmode\hbox{%
\textit{\,\,\,\,\, #1 Mathematical subject classification:\ }}}
\newcommand\blfootnote[1]{%
  \begingroup
  \renewcommand\thefootnote{}\footnote{#1}%
  \addtocounter{footnote}{-1}%
  \endgroup
}

\def \e {\varepsilon}

\def \phi {\varphi}

\def \R {\mathbb{R}}

\def \G{\Gamma}

\def \vf{\varphi}

\newcommand{\Rn}{\mathbb R^n}
\newcommand{\Rm}{\mathbb R^m}

\newcommand{\p}{\partial}

\newcommand{\la}{\lambda}

\numberwithin{equation}{section}

\newcommand{\beq}{\begin{equation}}
\newcommand{\bea}[1]{\begin{array}{#1} }
\newcommand{\eeq}{ \end{equation}}
\newcommand{\ea}{ \end{array}}

\newcommand{\sa}{\langle}
\newcommand{\da}{\rangle}

\newtheorem{theorem}{Theorem}[section]
\newtheorem{lemma}[theorem]{Lemma}

\newtheorem{remark}[theorem]{Remark}

\DeclareMathOperator{\sech}{sech}

\begin{document}

\title[]{Heat flow and the fundamental solution of the Baouendi-Grushin operator}

{\blfootnote{\MSC[2020]{Primary 35A08, 35K08; Secondary 35H20, 35J70, 33C10}}}


\date{}

\begin{abstract}
For the Baouendi-Grushin operator $\mathscr L_\gamma = \Delta_x + |x|^{2\gamma}\Delta_y$ in $\R^{m+k}$, with $\gamma>0$, we prove that the time integral of the heat kernel with pole at the origin equals $C_{\gamma,Q}\, N_\gamma^{2-Q}$, where $N_\gamma$ is the gauge and $Q$ the homogeneous dimension. When $\gamma\not=1$ no closed form of the heat kernel is known. After a partial Fourier transform in $y$, the time integral of the heat kernel becomes the Green function of an anharmonic oscillator, which solves a generalized modified Bessel equation, and the gauge appears when the inverse Fourier transform is computed as a Hankel transform. In particular, the heat flow provides a new proof that $C_{\gamma,Q}\, N_\gamma^{2-Q}$ is a fundamental solution of $-\mathscr L_\gamma$ with pole at the origin.
\end{abstract}

\author{Nicola Garofalo}
\address{School of Mathematical and Statistical Sciences\\ Arizona State University}\email[Nicola Garofalo]{nicola.garofalo@asu.edu}

\maketitle

\section{Introduction}\label{S:intro}

Consider the heat kernel $p(x,t) = (4\pi t)^{-\frac n2} e^{-\frac{|x|^2}{4t}}$ in $\Rn$, with $n\ge 3$, and denote by $\sigma_{n-1} = \frac{2\pi^{n/2}}{\G(n/2)}$ the $(n-1)$-dimensional measure of the unit sphere. A simple computation shows that for every $x\not= 0$
\[
\int_0^\infty p(x,t)\, dt = \frac{1}{(n-2)\sigma_{n-1}|x|^{n-2}},
\]
the fundamental solution of $-\Delta$ with pole at the origin. That the time integral of the heat kernel is a fundamental solution is a general principle, but here the explicit knowledge of the heat kernel also reveals the explicit form of the Newtonian potential. The same happens in the Heisenberg group, where the time integral of the heat kernel found by Gaveau \cite{Gav} and Hulanicki \cite{Hu} gives Folland's fundamental solution \cite{Fo}, and more generally in groups of Heisenberg type, where one obtains Kaplan's fundamental solution \cite{Ka}, see \cite{GTaim}. 
This note is about a similar conclusion for a class of degenerate differential operators whose heat kernels are, with one exception, not known in closed form.

Given $\gamma>0$ and $m, k\in \mathbb N$, we consider in $\Rm_x\times \R^k_y$ the partial differential operator
\begin{equation}\label{Pg}
\mathscr L_\gamma = \Delta_x + |x|^{2\gamma} \Delta_y,
\end{equation}
introduced by Baouendi \cite{Ba} and Grushin \cite{Gr}. For its role in analysis and geometry, and for some relevant literature, we refer to \cite{Gft, DGGS, BG}.
The operator \eqref{Pg} is translation invariant in $y$, and homogeneous of degree two with respect to the dilations $\delta_\rho(x,y) = (\rho x,\rho^{\gamma+1}y)$, for which $d(\delta_\rho(x,y)) = \rho^Q dx dy$, with
\begin{equation}\label{Q}
Q = m + (\gamma+1) k.
\end{equation}
Let $p_\gamma((x,y),t)$ denote the heat kernel for \eqref{Pg} with pole at $(0,0)$. When $\gamma\not=1$ no closed form of $p_\gamma$ is known, see \cite[pp.~189-190]{BFI}. 
The only exception is when $\gamma =1$. In \cite[Theorem 3.4]{GTpotan} the following heat kernel was found for the operator $\Delta_x + \frac{|x|^2}{4}\Delta_y$, see also the references therein:
\begin{align}\label{parBfs}
q((x,y),(x',y'),t) & =  \frac{2^k}{(4\pi t)^{\frac{m}2 +k}} \int_{\R^k} e^{- \frac it \langle \la,y'-y\rangle}   \left(\frac{|\la|}{\sinh |\la|}\right)^{\frac m2} 
\\
& \times  e^{-\frac{|\la|}{4t \tanh |\la|} ((|x|^2 +|x'|^2) - 2 \langle x,x'\rangle \sech |\la|)} d\la,
\notag\end{align}
and one has $P_1((x,y),(x',y'),t) = 2^{-k} q((x,\frac y2),(x',\frac{y'}2),t)$, with $P_1$ as in \eqref{PgK} below.
Notwithstanding the lack of a closed form of $p_\gamma$ when $\gamma\not=1$, our main result in this note is the following.

\medskip

\begin{theorem}\label{T:grushin}
Let $m, k\in \mathbb N$ and $\gamma>0$, and let $\nu = \frac{m-2}{2(\gamma+1)}$. For every $(x,y)\in \R^{m+k}$ with $x\not= 0$, one has
\begin{equation}\label{grushin}
\int_0^\infty p_\gamma((x,y),t)\, dt = C_{\gamma,Q}\, N_\gamma(x,y)^{2-Q},
\end{equation}
where $Q$ is as in \eqref{Q}, 
\begin{equation}\label{Ng}
N_\gamma(x,y) = \left(|x|^{2(\gamma+1)} + (\gamma+1)^2 |y|^2\right)^{\frac{1}{2(\gamma+1)}},
\end{equation}
and the constant in \eqref{grushin} is given by
\begin{equation}\label{const}
C_{\gamma,Q} = \frac{(\gamma+1)^{k-1}\, \G(\frac m2)\, \G(\frac{Q-2}{2(\gamma+1)})}{4\pi^{\frac{m+k}2}\, \G(\nu+1)}.
\end{equation}
\end{theorem}

\medskip

Concerning Theorem \ref{T:grushin} we mention that the gauge \eqref{Ng} was first introduced in the paper \cite{Gjde}. In Proposition 2.1 of the same paper it was proved that \begin{equation}\label{Gam}
\Gamma_\gamma(x,y) = C^\star_{\gamma,Q}\, N_\gamma(x,y)^{2-Q}
\end{equation}
is a fundamental solution of $-\mathscr L_\gamma$ with pole at the origin, where
\begin{equation}\label{CgQ}
[C^\star_{\gamma,Q}]^{-1} = (Q-2)(Q+2\gamma) \int_{\R^{m+k}} \frac{|x|^{2\gamma}}{\left[1 + N_\gamma(x,y)^{2(\gamma+1)}\right]^{1 + \frac{Q+2\gamma}{2(\gamma+1)}}}\ dx dy.
\end{equation}
The gauge \eqref{Ng} and the Baouendi-Grushin Newtonian potential \eqref{Gam} were instrumental in the discovery in \cite{Gjde} of a far-reaching monotonicity formula of Almgren type for the degenerate operators \eqref{Pg}.

Since in Section \ref{S:const} we show that $C^\star_{\gamma,Q} = C_{\gamma,Q}$, we note explicitly that Theorem \ref{T:grushin} says that the heat flow of $\mathscr L_\gamma$ discovers the gauge \eqref{Ng}. 
To provide the reader with an explanation of why this is possible, 
let $\mathscr L_\gamma$ also denote the Friedrichs extension of \eqref{Pg} from $C^\infty_0(\R^{m+k})$, and let $p_\gamma$ be the heat kernel of $\p_t - \mathscr L_\gamma$ with pole at the origin, defined in \eqref{pg} below through the partial Fourier transform in $y$ and the heat kernels $h_\la(x,x',t)$ of the anharmonic oscillators 
\[
H_\la = -\Delta_x + 4\pi^2|\la|^2|x|^{2\gamma},\ \ \ \  \la\in \R^k.
\]
When $\gamma = 1$ the function $h_\la$ is given by Mehler's formula, see \cite[Remark 3.3]{GTpotan}. When $\gamma\not=1$ no closed form of $h_\la$ is known, and, besides the general facts collected in Section \ref{S:hk}, the only information on $p_\gamma$ that we will use is the equation satisfied by its partial Fourier transform.

After such partial Fourier transform in $y$, the function $p_\gamma$ becomes the heat kernel $h_\la(x,0,t)$ of the anharmonic oscillator $H_\la$, and in the variable $s = \frac{|x|^{\gamma+1}}{\gamma+1}$ the operator $H_\la$ acting on functions which are radial in $x$ becomes, up to the factor $|x|^{2\gamma}$, a modified Bessel operator in the fractal dimension $\frac{m+2\gamma}{\gamma+1}$. Such factor cannot be removed from the evolution equation, in which its reciprocal multiplies $\p_t$, so that the evolution problem is not reduced to a Bessel heat equation. But it disappears from the stationary equation satisfied by $\int_0^\infty h_\la(x,0,t)\, dt$, the Green function of $H_\la$ with pole at the origin, which is therefore $s^{-\nu}$ times a Macdonald function of order $\nu$ in the variable $2\pi|\la| s$, see Section \ref{S:time}. The inverse Fourier transform in $\la$ is then a Hankel transform which can be computed in closed form, and it is at this stage that the gauge \eqref{Ng} appears, see Section \ref{S:freq}. The identity $C^\star_{\gamma,Q} = C_{\gamma,Q}$, which also follows from \cite[Lemma 3.4]{Gfinheat}, is obtained in Section \ref{S:const} by a direct evaluation of the integral in \eqref{CgQ}. Finally, by Lemma \ref{L:green} and Fubini's theorem, $\int_0^\infty p_\gamma\, dt$ is a fundamental solution of $-\mathscr L_\gamma$, and as we have mentioned above, the argument thus also gives an independent proof of \cite[Prop. 2.1]{Gjde}.

When $\gamma = 1$, Theorem \ref{T:grushin} is consistent with \cite[Corollary 1.3]{GTaim} and \cite[Theor. 1.3]{Gmehler}, see Remark \ref{R:consistent}. For the heat kernel of \eqref{Pg} when $\gamma = 1$ we refer to \cite{BFI11, CCFI, CCGGL, GTpotan, Ste}.
When $\gamma\in \mathbb N$, the operator \eqref{Pg} is a sum of squares of polynomial vector fields satisfying H\"ormander's condition, to which the lifting theorem in \cite{RoS} applies, and Bauer, Furutani and Iwasaki found in \cite[Cor. 7.4]{BFI} a fundamental solution of $\mathscr L_\gamma$ with pole at an arbitrary point, which for the pole at the origin reduces to \eqref{Gam}. Their approach integrates in time the explicit heat kernel of an auxiliary operator obtained after a partial Fourier transform in $y$ and polar coordinates in $x$, and uses $\gamma\in \mathbb N$ to sum a series of spherical harmonics; the heat kernel $p_\gamma$ does not appear in \cite{BFI}. See also \cite{BGG} for other special values of $\gamma$.

\medskip

\noindent \textbf{Notation}. Throughout the paper $m, k\in \mathbb N$ and $\gamma>0$ are fixed, $Q$ is as in \eqref{Q}, and we let
\begin{equation}\label{nutheta}
\nu = \frac{m-2}{2(\gamma+1)},\ \ \ \ \ \ \ \theta = \frac k2 + \nu = \frac{Q-2}{2(\gamma+1)}.
\end{equation}
We note that $\nu>-\frac 12$, that $\nu = 0$ if and only if $m=2$, and that $\theta>0$. We denote by $\omega_{m-1} = \frac{2\pi^{m/2}}{\G(m/2)}$ the surface measure of the unit sphere $\mathbb S^{m-1}\subset \Rm$, with the understanding that $\omega_0 = 2$. We will use the function
\begin{equation}\label{sdef}
s = s(x) = \frac{|x|^{\gamma+1}}{\gamma+1},
\end{equation}
in terms of which the gauge \eqref{Ng} can be written as
\begin{equation}\label{Ns}
s^2 + |y|^2 = \frac{N_\gamma(x,y)^{2(\gamma+1)}}{(\gamma+1)^2}.
\end{equation}
The partial Fourier transform with respect to $y\in \R^k$ is
\[
\hat u((x,\la),t) = \int_{\R^k} e^{-2\pi i\sa\la,y\da} u((x,y),t)dy.
\]
Finally, we say that a function $f$ on $\R^k$ is spherically symmetric if $f(\la) = f^\star(|\la|)$ for some $f^\star$; when $k = 1$ we simply mean that $f$ is even. We will use the following classical formula of Bochner, see \cite[Theor. 40 on p. 69]{BC}, which holds for $k\ge 1$ and for every spherically symmetric $f\in L^1(\R^k)$,
\begin{equation}\label{boch}
\hat f(\xi) = \frac{2\pi}{|\xi|^{\frac k2 - 1}} \int_0^\infty r^{\frac k2} f^\star(r)
J_{\frac{k}2-1}(2\pi |\xi| r) dr,\ \ \ \ \ \ \xi\not= 0.
\end{equation}
Since the inverse Fourier transform of a spherically symmetric function coincides with its Fourier transform, \eqref{boch} also computes $\check f$.

\section{The heat kernel with pole at the origin}\label{S:hk}

Since \eqref{Pg} is translation invariant in $y$, the Cauchy problem
\begin{equation}\label{cpg}
\p_t u - \mathscr L_\gamma u = 0,\ \ \ \ \ u((x,y),0) = \vf(x,y),
\end{equation}
is converted by the partial Fourier transform into the following one
\begin{equation}\label{cpgl}
\p_t \hat u + H_\la \hat u = 0,\ \ \ \ \ \hat u((x,\la),0) = \hat \vf(x,\la),
\end{equation}
where for $\la\in \R^k$ we have let
\begin{equation}\label{Hl}
H_\la = - \Delta_x + 4\pi^2 |\la|^2 |x|^{2\gamma}.
\end{equation}
We collect the facts about \eqref{Hl} that will be needed. Let $\la\not=0$. Since the potential $V(x) := 4\pi^2|\la|^2|x|^{2\gamma}$ is continuous, nonnegative, and tends to $\infty$ as $|x|\to \infty$, the operator $H_\la$ is self-adjoint and nonnegative on $L^2(\Rm)$, with compact resolvent, see \cite[Theor. XIII.67]{RS}. Its spectrum is thus discrete, and its first eigenvalue $\mu_1(\la)$ is strictly positive, since 
\[
\sa H_\la \psi,\psi\da = \int_{\Rm} |\nabla \psi|^2 + 4\pi^2 |\la|^2 \int_{\Rm} |x|^{2\gamma} \psi^2 = 0
\]
forces $\psi \equiv 0$. The change of variable $x\to |\la|^{-\frac{1}{\gamma+1}} x$ shows that $H_\la$ is unitarily equivalent to $|\la|^{\frac{2}{\gamma+1}} H_{\la_0}$, with $|\la_0| = 1$, and therefore
\begin{equation}\label{mu1}
\mu_1(\la) = \mu_1 |\la|^{\frac{2}{\gamma+1}},\ \ \ \ \ \ \mu_1 = \mu_1(\la_0)>0.
\end{equation}
The semigroup $e^{-tH_\la}$ admits a kernel $h_\la(x,x',t)$, continuous on $\Rm\times\Rm\times (0,\infty)$, symmetric in $(x,x')$, and such that
\begin{equation}\label{dom}
0\le h_\la(x,x',t) \le (4\pi t)^{-\frac m2} e^{-\frac{|x-x'|^2}{4t}},
\end{equation}
see \cite[Theor. 1.8.1 and Example 2.1.9]{Da}. For $\la = 0$ the kernel $h_0$ is the Gauss-Weierstrass kernel. The upper bound in \eqref{dom} implies in particular that $e^{-tH_\la}$ is a contraction in $L^\infty(\Rm)$. Since $H_\la$ commutes with the rotations of $\Rm$, one has for every $R\in \mathbb O(m)$,
\[
h_\la(Rx,Rx',t) = h_\la(x,x',t),
\]
and therefore the function $x\to h_\la(x,0,t)$ is spherically symmetric.

\begin{lemma}\label{L:decay}
Let $\la\in \R^k\setminus\{0\}$. For every $x, x'\in \Rm$ and $t>0$ one has
\[
h_\la(x,x',t) \le (2\pi t)^{-\frac m2} e^{-\frac{t}2 \mu_1 |\la|^{\frac{2}{\gamma+1}}}.
\]
\end{lemma}

\begin{proof}
By the semigroup property and the symmetry of the kernel, for $0<\tau<t$ we have $h_\la(x,\cdot,t) = e^{-(t-\tau)H_\la} h_\la(x,\cdot,\tau)$. Since $\|e^{-sH_\la}\|_{L^2\to L^2} = e^{-s\mu_1(\la)}$, we infer
\[
h_\la(x,x,2t) = ||h_\la(x,\cdot,t)||^2_{L^2(\Rm)} \le e^{-2(t-\tau)\mu_1(\la)} ||h_\la(x,\cdot,\tau)||^2_{L^2(\Rm)} = e^{-2(t-\tau)\mu_1(\la)} h_\la(x,x,2\tau).
\]
Taking $\tau = t/2$, and using \eqref{dom}, we find $h_\la(x,x,2t)\le e^{-t\mu_1(\la)} (4\pi t)^{-\frac m2}$. On the other hand, by Cauchy-Schwarz
\[
h_\la(x,x',2t) = \int_{\Rm} h_\la(x,\xi,t) h_\la(\xi,x',t) d\xi \le \left(h_\la(x,x,2t) h_\la(x',x',2t)\right)^{\frac 12} \le e^{-t\mu_1(\la)} (4\pi t)^{-\frac m2}.
\]
Replacing $t$ with $t/2$, and keeping \eqref{mu1} in mind, we reach the desired conclusion.

\end{proof}

Since $C^\infty_0(\Rm)$ is a form core for each $H_\la$, the closure of the form of $-\mathscr L_\gamma$ on $C^\infty_0(\R^{m+k})$ is, after the partial Fourier transform, the direct integral in $\la$ of the forms of $H_\la$, and therefore $\widehat{e^{t\mathscr L_\gamma}\vf}(\cdot,\la) = e^{-tH_\la}\hat\vf(\cdot,\la)$ for a.e. $\la$. Measurability in $\la$ follows from 
\[
h_\la(x,x',t) = |\la|^{\frac{m}{\gamma+1}} h_{\la_0}(|\la|^{\frac{1}{\gamma+1}}x,|\la|^{\frac{1}{\gamma+1}}x',|\la|^{\frac{2}{\gamma+1}}t).
\]
Hence the solution of \eqref{cpg} is given by
\begin{align*}
u((x,y),t) & = \int_{\R^k} e^{2\pi i \sa y,\la\da} \int_{\Rm} h_\la(x,x',t) \hat \vf(x',\la) dx' d\la
\\
& = \int_{\Rm}\int_{\R^k} P_\gamma((x,y),(x',y'),t) \vf(x',y') dy' dx',
\end{align*}
where
\begin{equation}\label{PgK}
P_\gamma((x,y),(x',y'),t) = \int_{\R^k} e^{2\pi i\sa y-y',\la\da} h_\la(x,x',t) d\la.
\end{equation}
By Lemma \ref{L:decay} the integral in \eqref{PgK} converges absolutely. Accordingly, the heat kernel of $\p_t - \mathscr L_\gamma$ with pole at the origin is the function
\begin{equation}\label{pg}
p_\gamma((x,y),t) = P_\gamma((x,y),(0,0),t) = \int_{\R^k} e^{2\pi i\sa y,\la\da} h_\la(x,0,t) d\la.
\end{equation}

\section{Integration in time at fixed frequency}\label{S:time}

For $\la\in \R^k\setminus\{0\}$ and $x\in \Rm\setminus\{0\}$ we define
\begin{equation}\label{w}
w(x,\la) = \int_0^\infty h_\la(x,0,t) dt.
\end{equation}
By \eqref{dom} and Lemma \ref{L:decay} we have
\begin{equation}\label{wbound}
0\le w(x,\la) \le \int_0^1 (4\pi t)^{-\frac m2} e^{-\frac{|x|^2}{4t}} dt + \int_1^\infty (2\pi t)^{-\frac m2} e^{-\frac t2 \mu_1 |\la|^{\frac{2}{\gamma+1}}} dt,
\end{equation}
which shows that the integral in \eqref{w}:
\begin{itemize}
\item[(a)] converges for $x\not= 0$,
\item[(b)] that $w(\cdot,\la)\in L^1_{loc}(\Rm)$ (the first integral in the right-hand side of \eqref{wbound} is bounded by $C |x|^{2-m}$ when $m\ge 3$,  by $C(1+|\log |x||)$ when $m = 2$, and  by a constant when $m=1$),
\item[(c)] and that, by dominated convergence (for each $t>0$ one has $h_\la(x,0,t)\to 0$ as $|x|\to\infty$ by \eqref{dom}, and for $|x|\ge 1$ the integrand is dominated by the minimum of the two bounds in \eqref{dom} and Lemma \ref{L:decay}, which is integrable on $(0,\infty)$),
\begin{equation}\label{winfty}
w(x,\la) \to 0\ \ \ \ \text{as}\ |x|\to \infty.
\end{equation}
\end{itemize}
Moreover, $w(\cdot,\la)$ is spherically symmetric, since so is $h_\la(\cdot,0,t)$.

\begin{lemma}\label{L:green}
Let $\la\in \R^k\setminus\{0\}$. For every $\psi\in C^\infty_0(\Rm)$ one has
\begin{equation}\label{green}
\int_{\Rm} w(x,\la) H_\la \psi(x) dx = \psi(0).
\end{equation}
\end{lemma}

\begin{proof}
Since $\psi$ belongs to the domain of $H_\la$, we have $\p_t e^{-tH_\la}\psi = - e^{-tH_\la} H_\la \psi$. Using the symmetry of the kernel $h_\la$, for $0<\e<T<\infty$ we thus find
\[
\int_\e^T \int_{\Rm} h_\la(x,0,t) H_\la \psi(x) dx dt = \int_\e^T \left(e^{-tH_\la} H_\la \psi\right)(0) dt = \left(e^{-\e H_\la} \psi\right)(0) - \left(e^{-T H_\la} \psi\right)(0).
\]
By Lemma \ref{L:decay} we have 
\[
|(e^{-T H_\la} \psi)(0)|\le (2\pi T)^{-\frac m2} e^{-\frac T2 \mu_1(\la)} ||\psi||_{L^1(\Rm)}\to 0,\ \ \ \text{as}\ T\to \infty.
\]
 On the other hand, since $e^{-sH_\la}$ is a contraction in $L^\infty(\Rm)$, we have
\[
\left|\left(e^{-\e H_\la} \psi\right)(0) - \psi(0)\right| = \left|\int_0^\e \left(e^{-sH_\la} H_\la \psi\right)(0) ds\right| \le \e\ ||H_\la \psi||_{L^\infty(\Rm)} \to 0\ \ \ \text{as}\ \e\to 0^+.
\]
Finally, since $h_\la\ge 0$ and $w(\cdot,\la)\in L^1_{loc}(\Rm)$, Tonelli's theorem gives
\[
\int_0^\infty\int_{\Rm} h_\la(x,0,t) |H_\la \psi(x)| dx dt \le ||H_\la \psi||_{L^\infty(\Rm)} \int_{\text{supp}\ \psi} w(x,\la) dx<\infty,
\]
and we can let $\e\to 0^+$ and $T\to \infty$ in the above identity by dominated convergence. This proves \eqref{green}.

\end{proof}

We now solve \eqref{green}. Lemma \ref{L:green} says that $w(\cdot,\la)$ is a distributional solution of $H_\la w = \delta$ in $\Rm$, where $\delta$ is the Dirac delta at the origin. In particular, since the potential in \eqref{Hl} is $C^\infty$ in $\Rm\setminus\{0\}$, we have $H_\la w = 0$ there, and by hypoellipticity $w(\cdot,\la)\in C^\infty(\Rm\setminus\{0\})$. By the spherical symmetry of $w(\cdot,\la)$, we can write $w(x,\la) = W(|x|)$, where $W\in C^\infty(0,\infty)$ solves the ordinary differential equation
\begin{equation}\label{odeW}
W'' + \frac{m-1}{r} W' - 4\pi^2 |\la|^2 r^{2\gamma} W = 0,\ \ \ \ \ \ \ r>0.
\end{equation}
If we make the change of variables
\begin{equation}\label{cv}
W(r) = r^{1-\frac m2} \Phi(\zeta),\ \ \ \ \ \ \zeta = \frac{2\pi |\la|}{\gamma+1} r^{\gamma+1} = 2\pi |\la| s,
\end{equation}
a direct computation shows that
\[
\frac{r^{\frac m2 + 1}}{(\gamma+1)^2}\left(W'' + \frac{m-1}{r} W' - 4\pi^2 |\la|^2 r^{2\gamma} W\right) = \zeta^2 \Phi''(\zeta) + \zeta \Phi'(\zeta) - (\zeta^2 + \nu^2) \Phi(\zeta),
\]
with $\nu$ as in \eqref{nutheta}. Therefore, \eqref{odeW} is equivalent to the modified Bessel equation of order $\nu$ for $\Phi$, whose general solution is $\Phi(\zeta) = A K_\nu(\zeta) + B I_\nu(\zeta)$, see \cite[Sec. 5.7]{Le}. Here, $I_\nu$ and $K_\nu$ respectively denote the modified Bessel functions of the first and third kind, and we recall that $K_{-\nu} = K_\nu$. Since as $\zeta\to \infty$ one has $I_\nu(\zeta)\sim (2\pi\zeta)^{-\frac12} e^{\zeta}$, whereas $K_\nu(\zeta)\sim (\frac{\pi}{2\zeta})^{\frac12} e^{-\zeta}$, \eqref{winfty} forces $B = 0$. We conclude that
\begin{equation}\label{wK}
w(x,\la) = A(\la)\ |x|^{1-\frac m2} K_\nu(2\pi |\la| s),
\end{equation}
for some constant $A(\la)$ which we now determine using \eqref{green}. Fix $\psi\in C^\infty_0(\Rm)$ with $\psi(0) = 1$. Since $H_\la w = 0$ in $\Rm\setminus\{0\}$, integrating by parts in $\{|x|>\e\}$ we find
\begin{align*}
\int_{\{|x|>\e\}} w H_\la \psi dx & = \int_{\{|x|>\e\}} (\psi \Delta w - w \Delta \psi) dx = \int_{\{|x| = \e\}} \left(w\, \p_r \psi - \psi\, \p_r w \right) dS_\e
\\
& = \e^{m-1} W(\e) \int_{\mathbb S^{m-1}} \p_r \psi(\e \omega) d\omega - \e^{m-1} W'(\e) \int_{\mathbb S^{m-1}} \psi(\e \omega) d\omega.
\end{align*}
When $m = 1$ the integrals over $\mathbb S^0 = \{\pm 1\}$ are to be interpreted as sums. Letting $\e\to 0^+$ in the latter identity, and using \eqref{green} and $w\in L^1_{loc}(\Rm)$, we obtain
\begin{equation}\label{flux}
\lim_{\e\to 0^+} \e^{m-1} W(\e) \int_{\mathbb S^{m-1}} \p_r \psi(\e \omega) d\omega - \omega_{m-1} \lim_{\e\to 0^+} \e^{m-1} W'(\e)  = 1,
\end{equation}
provided that the two limits exist. We now use the behavior of $K_\nu$ near $\zeta = 0$. We treat separately three cases.

\medskip

\noindent \emph{Case $m\ge 3$}. In this case $\nu>0$, and we have $K_\nu(\zeta) = \frac 12 \G(\nu) (\zeta/2)^{-\nu} (1 + o(1))$ and $K'_\nu(\zeta) = - \frac{\nu}{2} \G(\nu) (\zeta/2)^{-\nu} \zeta^{-1} (1 + o(1))$ as $\zeta\to 0^+$, see \cite[Sec. 5.7]{Le}. Since by \eqref{nutheta} we have $(\gamma+1)\nu = \frac{m-2}{2}$, from \eqref{cv} we find
\begin{align}\label{Wasy}
W(r) & = \frac{A(\la)}{2} \G(\nu) \left(\frac{\pi|\la|}{\gamma+1}\right)^{-\nu} r^{2-m} \left(1+o(1)\right),
\\
W'(r) & = (2-m) \frac{A(\la)}{2} \G(\nu) \left(\frac{\pi|\la|}{\gamma+1}\right)^{-\nu} r^{1-m} \left(1+o(1)\right),
\notag
\end{align}
as $r\to 0^+$. The first limit in \eqref{flux} vanishes, and we obtain
\begin{equation}\label{Am}
A(\la) = \frac{2}{(m-2)\, \omega_{m-1}\, \G(\nu)} \left(\frac{\pi|\la|}{\gamma+1}\right)^{\nu}.
\end{equation}

\medskip

\noindent \emph{Case $m = 2$}. In this case $\nu = 0$, and $K_0(\zeta) = - \log \zeta + O(1)$, $K_0'(\zeta) = - K_1(\zeta) = - \zeta^{-1} (1 + o(1))$ as $\zeta\to 0^+$, see \cite[Sec. 5.7]{Le}. Since $\zeta = \frac{2\pi|\la|}{\gamma+1} r^{\gamma+1}$, from \eqref{cv} we obtain $W(r) = - (\gamma+1) A(\la) \log r + O(1)$, and $W'(r) = - (\gamma+1) A(\la) r^{-1} (1 + o(1))$ as $r\to 0^+$. Again, the first limit in \eqref{flux} vanishes, and we find
\begin{equation}\label{A2}
A(\la) = \frac{1}{(\gamma+1)\omega_1} = \frac{1}{2\pi(\gamma+1)}.
\end{equation}

\medskip

\noindent \emph{Case $m = 1$}. In this case $-\frac 12<\nu = - \frac{1}{2(\gamma+1)}<0$, and $K_\nu = K_{|\nu|}$, with
\[
K_{|\nu|}(\zeta) = \frac 12 \G(|\nu|) \left(\frac{\zeta}{2}\right)^{-|\nu|} \left(1 + O(\zeta^2)\right) + \frac 12 \G(-|\nu|) \left(\frac{\zeta}{2}\right)^{|\nu|}\left(1 + O(\zeta^2)\right),
\]
as $\zeta\to 0^+$, which follows from the series expansion of $I_{\pm|\nu|}$, from $K_{|\nu|} = \frac{\pi}{2} \frac{I_{-|\nu|} - I_{|\nu|}}{\sin \pi |\nu|}$, and from the reflection formula for the gamma function, see \cite[Sec. 5.7]{Le}. Since now $(\gamma+1)|\nu| = \frac 12$, from \eqref{cv} we find
\[
W(r) = \frac{A(\la)}2 \G(|\nu|) \left(\frac{\pi|\la|}{\gamma+1}\right)^{-|\nu|} + \frac{A(\la)}2 \G(-|\nu|) \left(\frac{\pi|\la|}{\gamma+1}\right)^{|\nu|} r + O(r^{2\gamma+2}),
\]
which shows that $W$ is continuous up to $r = 0$, and that $W'(0^+) = \frac{A(\la)}2 \G(-|\nu|) (\frac{\pi|\la|}{\gamma+1})^{|\nu|}$. Since $\p_r \psi(\e) + \p_r\psi(-\e) = \psi'(\e) - \psi'(-\e)\to 0$, the first limit in \eqref{flux} vanishes also in this case, and we obtain $A(\la) \G(-|\nu|) (\frac{\pi|\la|}{\gamma+1})^{|\nu|} = -1$. Recalling that $\nu = -|\nu|$, we have thus found
\begin{equation}\label{A1}
A(\la) = - \frac{1}{\G(\nu)} \left(\frac{\pi|\la|}{\gamma+1}\right)^{\nu}.
\end{equation}

\medskip

Since $\omega_0 = 2$, we see that \eqref{A1} is \eqref{Am} with $m = 1$. Furthermore, using $\G(\nu+1) = \nu \G(\nu)$, and $\frac{\nu}{\frac m2 - 1} = \frac{1}{\gamma+1}$, we can rewrite \eqref{Am} in the following form, which for $m = 2$ returns \eqref{A2},
\begin{equation}\label{Afinal}
A(\la) = \frac{\G(\frac m2)}{2(\gamma+1) \pi^{\frac m2}\, \G(\nu+1)} \left(\frac{\pi}{\gamma+1}\right)^{\nu} |\la|^{\nu}.
\end{equation}
Substituting \eqref{Afinal} in \eqref{wK}, we conclude that for every $m\in \mathbb N$ one has
\begin{equation}\label{wfinal}
w(x,\la) = \frac{\G(\frac m2)}{2(\gamma+1) \pi^{\frac m2}\, \G(\nu+1)} \left(\frac{\pi}{\gamma+1}\right)^{\nu} |\la|^{\nu}\ |x|^{1-\frac m2}\ K_\nu\left(\frac{2\pi |\la| |x|^{\gamma+1}}{\gamma+1}\right).
\end{equation}

\section{Integration in the frequency variable}\label{S:freq}

We now fix $x\in \Rm\setminus\{0\}$ and compute the inverse Fourier transform of \eqref{wfinal} with respect to $\la$. We first observe that $w(x,\cdot)\in L^1(\R^k)$. The exponential decay of $K_\nu$ takes care of the integrability at infinity. Near $\la = 0$, from \eqref{wfinal} and the expansions of $K_\nu$ recalled above, one has that $w(x,\la)$ is bounded when $m\ge 3$, that $|w(x,\la)|\le C(1+|\log|\la||)$ when $m = 2$, and that $|w(x,\la)|\le C |\la|^{2\nu} = C|\la|^{-\frac{1}{\gamma+1}}$ when $m = 1$. In all cases, $w(x,\cdot)$ is integrable near $\la = 0$ since $k\ge 1>\frac{1}{\gamma+1}$.

Since $h_\la\ge 0$, we have from \eqref{pg}, \eqref{w} and Tonelli's theorem
\[
\int_0^\infty \int_{\R^k} \left|e^{2\pi i\sa y,\la\da} h_\la(x,0,t)\right| d\la dt = \int_{\R^k} w(x,\la) d\la<\infty.
\]
By Fubini's theorem we thus obtain
\begin{equation}\label{Efourier}
\mathscr E_\gamma(x,y) = \int_0^\infty p_\gamma((x,y),t) dt = \int_{\R^k} e^{2\pi i\sa y,\la\da} w(x,\la) d\la.
\end{equation}
In particular, $\mathscr E_\gamma(x,\cdot)$ is continuous in $\R^k$. Since $\Gamma_\gamma$ is also continuous in $y$ when $x\not=0$, it suffices to consider $y\not= 0$. Since $w(x,\cdot)$ is spherically symmetric in $\R^k$, the inverse Fourier transform in \eqref{Efourier} is computed by \eqref{boch}. Writing $w(x,\la) = A_0 |\la|^\nu |x|^{1-\frac m2} K_\nu(2\pi |\la| s)$, with $s$ as in \eqref{sdef} and
\begin{equation}\label{A0}
A_0 = \frac{\G(\frac m2)}{2(\gamma+1) \pi^{\frac m2}\, \G(\nu+1)} \left(\frac{\pi}{\gamma+1}\right)^{\nu},
\end{equation}
we obtain from \eqref{Efourier}
\begin{equation}\label{Ehankel}
\mathscr E_\gamma(x,y) = \frac{2\pi A_0 |x|^{1-\frac m2}}{|y|^{\frac k2 - 1}} \int_0^\infty \rho^{\frac k2 + \nu} K_\nu(2\pi s \rho) J_{\frac k2 - 1}(2\pi |y| \rho) d\rho.
\end{equation}
At this point we recall the following classical formula, see \cite[Sec. 13.45]{Wa}, or also \cite[6.576 (3)]{GR},
\begin{align}\label{WaK}
\int_0^\infty \tau^{-\eta} K_\mu(a \tau) J_{\nu'}(b \tau) d\tau & = \frac{b^{\nu'} \G(\frac{\nu'-\eta+\mu+1}{2})\G(\frac{\nu'-\eta-\mu+1}{2})}{2^{\eta+1} a^{\nu'-\eta+1} \G(\nu'+1)} 
\\
& \times F\left(\frac{\nu'-\eta+\mu+1}{2}, \frac{\nu'-\eta-\mu+1}{2}; \nu'+1; - \frac{b^2}{a^2}\right),
\notag
\end{align}
valid for $a>0$, $b\ge 0$, and $\nu'-\eta+1>|\mu|$. We apply \eqref{WaK} with the choice
\[
\mu = \nu,\ \ \ \ \nu' = \frac k2 - 1,\ \ \ \ \eta = - \frac k2 - \nu,\ \ \ \ a = 2\pi s,\ \ \ \ b = 2\pi |y|.
\]
The condition $\nu'-\eta+1 = k + \nu >|\nu|$ is satisfied, since $k\ge 1>\frac{1}{\gamma+1} = 2|\nu|$ when $m = 1$, and $\nu\ge 0$ when $m\ge 2$. Keeping \eqref{nutheta} in mind, we find
\[
\frac{\nu'-\eta+\mu+1}{2} = \frac k2 + \nu = \theta,\ \ \ \ \ \ \ \frac{\nu'-\eta-\mu+1}{2} = \frac k2 = \nu' + 1.
\]
The hypergeometric function in \eqref{WaK} is therefore in the special form $F(a,b;b;x) = (1-x)^{-a}$, see \cite[(4) on p. 101]{E}, and we have
\[
F\left(\theta, \frac k2; \frac k2; - \frac{|y|^2}{s^2}\right) = \left(1 + \frac{|y|^2}{s^2}\right)^{-\theta} = s^{2\theta} (s^2 + |y|^2)^{-\theta}.
\]
We thus obtain from \eqref{WaK}
\begin{equation}\label{hankelvalue}
\int_0^\infty \rho^{\frac k2 + \nu} K_\nu(2\pi s \rho) J_{\frac k2 - 1}(2\pi |y| \rho) d\rho = 2^{\frac k2 + \nu - 1}  (2\pi |y|)^{\frac k2 - 1} (2\pi s)^{-k-\nu}\, \G(\theta)\, s^{2\theta} (s^2 + |y|^2)^{-\theta}.
\end{equation}
Substituting \eqref{hankelvalue} in \eqref{Ehankel}, and observing that $2\theta - k - \nu = \nu$, we see that the powers of $|y|$ cancel, and we are left with
\begin{equation}\label{E6}
\mathscr E_\gamma(x,y) = \frac{A_0\, \G(\theta)}{2\pi^{\frac k2 + \nu}}\ |x|^{1-\frac m2} s^{\nu}\ (s^2 + |y|^2)^{-\theta}.
\end{equation}
Now, by \eqref{sdef} and $(\gamma+1)\nu = \frac{m-2}2$, we have $|x|^{1-\frac m2} s^{\nu} = (\gamma+1)^{-\nu}$, and by \eqref{Ns} and $2(\gamma+1)\theta = Q-2$, we have $(s^2 + |y|^2)^{-\theta} = (\gamma+1)^{2\theta} N_\gamma(x,y)^{2-Q}$. Substituting these identities and \eqref{A0} in \eqref{E6}, and using $2\theta - 2\nu = k$, we finally obtain
\begin{equation}\label{Egamma}
\mathscr E_\gamma(x,y) = \frac{(\gamma+1)^{k-1}\, \G(\frac m2)\, \G(\theta)}{4\pi^{\frac{m+k}2}\, \G(\nu+1)}\ N_\gamma(x,y)^{2-Q}.
\end{equation}
Since $\theta = \frac{Q-2}{2(\gamma+1)}$, the constant in the right-hand side of \eqref{Egamma} is $C_{\gamma,Q}$ in \eqref{const}, and this proves \eqref{grushin}. In the next section we show that $C_{\gamma,Q} = C^\star_{\gamma,Q}$.

\section{Computing the constant}\label{S:const}

In this section we show that $C^\star_{\gamma,Q} = C_{\gamma,Q}$, where $C^\star_{\gamma,Q}$ is defined by \eqref{CgQ}. This also follows from \cite[Lemma 3.4]{Gfinheat} with $p=2$ and Euclidean norms, after the change of variable $y\to y/2$; for completeness we give the direct evaluation of the integral in \eqref{CgQ}. By \eqref{nutheta} we have
\[
1 + \frac{Q+2\gamma}{2(\gamma+1)} = \theta + 2,\ \ \ \ \ \ \ \ \ \theta + 2 - \frac k2 = \nu + 2.
\]
Using the elementary formula $\int_{\R^k} (a + |\xi|^2)^{-e} d\xi = \pi^{\frac k2} \frac{\G(e - \frac k2)}{\G(e)} a^{\frac k2 - e}$, valid for $a>0$ and $e>\frac k2$, with $a = 1 + |x|^{2(\gamma+1)}$ and $e = \theta + 2$, we first integrate in $y$ in \eqref{CgQ}, obtaining
\[
\int_{\R^{k}} \frac{dy}{\left[1 + |x|^{2(\gamma+1)} + (\gamma+1)^2 |y|^2\right]^{\theta+2}} = \frac{\pi^{\frac k2}}{(\gamma+1)^{k}}\frac{\G(\nu+2)}{\G(\theta+2)}\left(1 + |x|^{2(\gamma+1)}\right)^{-\nu-2}.
\]
Next, we integrate in $x$. The change of variable $u = r^{2(\gamma+1)}$, and the identity $\frac{m + 2\gamma}{2(\gamma+1)} - 1 = \nu$, give
\begin{align*}
\int_{\Rm} \frac{|x|^{2\gamma}\, dx}{\left(1 + |x|^{2(\gamma+1)}\right)^{\nu+2}} & = \omega_{m-1} \int_0^\infty \frac{r^{m-1+2\gamma}\, dr}{\left(1 + r^{2(\gamma+1)}\right)^{\nu+2}}
\\
& = \frac{\omega_{m-1}}{2(\gamma+1)} \int_0^\infty \frac{u^{\nu}\, du}{(1+u)^{\nu+2}} = \frac{\omega_{m-1}}{2(\gamma+1)(\nu+1)},
\end{align*}
where in the last equality we have used $\int_0^\infty u^{\nu}(1+u)^{-\nu-2} du = B(\nu+1,1) = \frac{1}{\nu+1}$. Combining the last two formulas, we find
\[
[C^\star_{\gamma,Q}]^{-1} = (Q-2)(Q+2\gamma)\ \frac{\pi^{\frac k2}\, \omega_{m-1}}{2(\gamma+1)^{k+1}}\ \frac{\G(\nu+1)}{\G(\theta+2)}.
\]
Since $(Q-2)(Q+2\gamma) = 4(\gamma+1)^2 \theta(\theta+1)$ and $\G(\theta+2) = \theta(\theta+1)\G(\theta)$, we conclude
\[
[C^\star_{\gamma,Q}]^{-1} = \frac{2\pi^{\frac k2}\, \omega_{m-1}\, \G(\nu+1)}{(\gamma+1)^{k-1}\, \G(\theta)} = \frac{4\pi^{\frac{m+k}2}\, \G(\nu+1)}{(\gamma+1)^{k-1}\, \G(\frac m2)\, \G(\theta)},
\]
which is the reciprocal of $C_{\gamma,Q}$ in \eqref{const}. This proves $C^\star_{\gamma,Q} = C_{\gamma,Q}$. Finally, by Lemma \ref{L:green} applied to $\psi=\hat\vf(\cdot,\la)$ and Fubini's theorem, $\mathscr E_\gamma$ is a fundamental solution of $-\mathscr L_\gamma$ with pole at the origin, and together with \eqref{Egamma} this gives an independent proof of \cite[Prop. 2.1]{Gjde}.

\medskip

\begin{remark}\label{R:consistent}
When $\gamma = 1$, Mehler's formula for $h_\la(x,0,t)$, see \cite[Remark 3.3]{GTpotan}, and the change of variable $\la\to \frac{\la}{4\pi t}$ in \eqref{pg} give $2^k p_1((x,2y),t) = G_{\frac m2,\frac m2+k}((x,y),t)$, where $G_{\alpha,\beta}$ is the generalized Mehler kernel in \cite[Prop. 1.2]{Gmehler}. With $m = 2\alpha$ one has $N_1(x,2y) = (|x|^4+16|y|^2)^{1/4}$, $Q = 2\alpha + 2k$, $\nu = \frac{\alpha-1}{2}$ and $\theta = \frac{\alpha+k-1}{2}$, and by Legendre's duplication formula $2^k$ times the constant in \eqref{const} equals the constant in \cite[Theor. 1.3]{Gmehler}, and in \cite[Corollary 1.3]{GTaim} when $m$ is the dimension of the horizontal layer of a group of Heisenberg type.
\end{remark}

\vskip 0.2in

\section*{Declarations}

\smallskip
\noindent\textbf{Competing interests.} The author has no competing
interests to declare that are relevant to the content of this article.

\smallskip

\noindent \textbf{Use of AI:} The author acknowledges the assistance of AI in organizing and revising the manuscript, and in reviewing the proofs and the language. The research ideas and the mathematical arguments are exclusively the author's, and he takes full responsibility for the content and the final text.


\bibliographystyle{amsplain}

\begin{thebibliography}{10}

\bibitem{Ba}
M. S. Baouendi, \emph{Sur une classe d'op\'erateurs elliptiques d\'eg\'en\'er\'es}. Bull. Soc. Math. France 95 (1967), 45-87.

\bibitem{BG}
A. Banerjee \& N. Garofalo, \emph{An improved Rellich type estimate for the eigenfunctions of the Baouendi-Grushin operator}, to appear in J. Anal. Math. 

\bibitem{BFI11}
W. Bauer, K. Furutani \& C. Iwasaki, \emph{Spectral analysis and geometry of sub-Laplacian and related Grushin-type operators}. Partial differential equations and spectral theory, 183-290, Oper. Theory Adv. Appl., 211, Birkh\"auser/Springer Basel AG, Basel, 2011.

\bibitem{BFI}
W. Bauer, K. Furutani \& C. Iwasaki, \emph{Fundamental solution of a higher step Grushin type operator}. Adv. Math. 271 (2015), 188-234.

\bibitem{BGG}
R. Beals, B. Gaveau \& P. Greiner, \emph{Green's functions for some highly degenerate elliptic operators}. J. Funct. Anal. 165 (1999), no. 2, 407-429.

\bibitem{BC}
S. Bochner \& K. Chandrasekharan, \emph{Fourier Transforms}, Annals of Mathematics Studies, no. 19, Princeton University Press, Princeton, N. J.; Oxford University Press, London, 1949. ix+219 pp.



\bibitem{CCFI}
O. Calin, D.-C. Chang, K. Furutani \& C. Iwasaki, \emph{Heat kernels for elliptic and sub-elliptic operators. Methods and techniques}. Applied and Numerical Harmonic Analysis. Birkh\"auser/Springer, New York, 2011. xviii+433 pp.

\bibitem{CCGGL}
C.-H. Chang, D.-C. Chang, B. Gaveau, P. Greiner \& H.-P. Lee, \emph{Geometric analysis on a step 2 Grusin operator}. Bull. Inst. Math. Acad. Sin. (N.S.) 4 (2009), no. 2, 119-188.

\bibitem{Da}
E. B. Davies, \emph{Heat kernels and spectral theory}. Cambridge Tracts in Mathematics, 92. Cambridge University Press, Cambridge, 1989. x+197 pp.

\bibitem{DGGS}
F. Dragoni, N. Garofalo, G. Giovannardi \& P. Salani, \emph{A fundamental solution for a subelliptic operator in Finsler geometry}. Indiana Univ. Math. J. 75 (2026), no. 1, 201-226.



\bibitem{E}
A. Erd\'elyi, W. Magnus, F. Oberhettinger \& F. G. Tricomi, \emph{Higher transcendental functions}. Vol. I. Based on notes left by Harry Bateman. With a preface by Mina Rees. With a foreword by E. C. Watson. Reprint of the 1953 original. Robert E. Krieger Publishing Co., Inc., Melbourne, Fla., 1981. xiii+302 pp.

\bibitem{Fo}
G. B. Folland, \emph{A fundamental solution for a subelliptic operator}. Bull. Amer. Math. Soc. 79~(1973), 373-376.


\bibitem{Gjde}
N. Garofalo, \emph{Unique continuation for a class of elliptic operators which degenerate on a manifold of arbitrary codimension}. J. Differential Equations 104 (1993), no. 1, 117-146.

\bibitem{Gft}
N. Garofalo, \emph{Fractional thoughts}.
Contemp. Math., 723
American Mathematical Society, [Providence], RI, 2019, 1-135.

\bibitem{Gfinheat}
N. Garofalo, \emph{On an evolution equation in sub-Finsler geometry}. Anal. Geom. Metr. Spaces 12 (2024), Paper No. 20240012.





\bibitem{Gmehler}
N. Garofalo, \emph{Fractal Mehler kernels and nonlinear geometric flows}. Preprint, arXiv:2512.23830. 

\bibitem{GTaim}
N. Garofalo \& G. Tralli, \emph{Feeling the heat in a group of Heisenberg type}.
Adv. Math. 381 (2021), Paper No. 107635, 42 pp.

\bibitem{GTpotan}
N. Garofalo \& G. Tralli, \emph{ Heat kernels for a class of hybrid evolution equations}.
Potential Anal. 59 (2023), no. 2, 823-856.

\bibitem{Gav}
B. Gaveau, \emph{Principe de moindre action, propagation de la chaleur et estim\'ees sous elliptiques sur certains groupes nilpotents}. Acta Math. 139~(1977), no. 1-2, 95-153.

\bibitem{GR}
I. S. Gradshteyn \& I. M. Ryzhik, \emph{Table of integrals, series, and products}. Seventh edition. Elsevier/Academic Press, Amsterdam, 2007. xlviii+1171 pp.



\bibitem{Gr}
V. V. Grushin, \emph{A certain class of hypoelliptic operators}. (Russian) Mat. Sb. (N.S.) 83 (125) (1970), 456-473.

\bibitem{Hu}
A. Hulanicki, \emph{The distribution of energy in the Brownian motion in the Gaussian field and analytic-hypoellipticity of certain subelliptic operators on the Heisenberg group}. Studia Math. 56~(1976), no. 2, 165-173.



\bibitem{Ka}
A. Kaplan, \emph{Fundamental solutions for a class of hypoelliptic PDE generated by composition of quadratic forms}. Trans. Amer. Math. Soc. 258~(1980), no. 1, 147-153.

\bibitem{Le}
N. N. Lebedev, \emph{Special functions and their applications}. Revised edition, translated from the Russian and edited by R. A. Silverman. Unabridged and corrected republication. Dover Publications, Inc., New York, 1972. xii+308 pp.



\bibitem{RS}
M. Reed \& B. Simon, \emph{Methods of modern mathematical physics. IV. Analysis of operators}. Academic Press, New York-London, 1978. xv+396 pp.

\bibitem{RoS}
L. P. Rothschild \& E. M. Stein, \emph{Hypoelliptic differential operators and nilpotent groups}. Acta Math. 137~(1976), no. 3-4, 247-320.



\bibitem{Ste}
K. Stempak, \emph{A transform for the Grushin operator with applications}.
J. Differential Equations 481 (2026), Paper No. 114630, 23 pp.

\bibitem{Wa}
G. N. Watson, \emph{A treatise on the theory of Bessel functions}. Reprint of the second (1944) edition. Cambridge Mathematical Library. Cambridge University Press, Cambridge, 1995. viii+804 pp.

\end{thebibliography}

\end{document}